\documentstyle[amscd,amssymb,verbatim]{amsart}

\theoremstyle{plain}
\newtheorem{Thm}{Theorem}
\newtheorem{Prop}{Proposition}[section]
\newtheorem{Lem}[Prop]{Lemma}

\newtheorem{Cor}{Corollary}

\theoremstyle{definition}

\begin{document}
\title[Integral of the exponential Brownian motion]
{On the distribution of the integral of the exponential Brownian motion }
\author{Leonid Tolmatz}
\address{Department of Mathematical Sciences, Xi'an Jiaotong--Liverpool University \\
  111 Ren Ai Road, Dushu Lake Higher Education Town, Suzhou Industrial Park, Suzhou, Jiangsu  
 215123, P.R. China}
\email{leonid.tolmatz@@xjtlu.edu.cn}

\date{June 23, 2008}


\begin{abstract}

The density distribution function of the integral of the exponential Brownian motion is determined explicitly in the form of a rapidly convergent series.

\end{abstract}

\maketitle

\section{Introduction} \label{S:intro}

Let \(B(s)\), \( 0 \le  s\le t \) denote a standard Brownian 
motion. 

In the present paper we consider the functional

\begin{equation}\label{E:fun}
\int_0^t e^{2 B(s)} ds.  
\end{equation}

This functional naturally appears in mathematical finance in connection with Asian options and in various other applications. 
Understanding of a functional is stipulated by computability of its distribution, and this problem  has attracted attention of researchers since 1992 when its study was initiated by M.Yor.
For details and further references see Matsumoto and Yor (2005a), (2005b), Dufresne (2000), Linetsky (2004), Schr\"oder (2006). \\

In the present article 
we obtain a new representation 
of the density distribution function $f(\lambda, t)$ of this functionl in the form of a rapidly convergent series. This series is given in Theorem \ref{Th:main} Section \ref{S:lt}. It allows computations with any desired accuracy in the whole domain, except small 
values of the argument, but for these asymptotic expansions are available; these will be discussed elsewhere. \\
These results were obtained by the author  between 1997 and 1999, but remained unpublished. In 2003 they were announced on the author's web page with some additional results and were communicated at the AARMS workshop on Financial Mathematics in St.Johns in 2003. \\
The method  applied readily leads to similar series for the cumulative distribution function and pricing of Asian options.
Our approach is based on the classical method of M.Kac. \\[10pt]

{\bf Notation and Conventions.} \\
$\zeta(z)$ is Riemann's Zeta function, \\
${}_1F_1(\alpha; \gamma;z) $ is the confluent hypergeometric function, \\
$H_m(x) $ are Hermite polynomials, \\
$ {\bf C} = 0.57721 \ldots $ is Euler's constant, \\
$[x]$ is the integer part of $x$, \\
$F(\lambda,t)$ and $f(\lambda,t)$ are the distribution and density distribution functions 
of functional (\ref{E:fun}).

\section{Double Laplace Transforms and the Corresponding Differential Equation.} \label{S:DLT}
Let

\begin{equation} \label{E:firstDist}
	F(\lambda,t) := P    \bigg\{\int_{0}^{t} e^{2 B(s)} ds <\lambda  \, \bigg\}
\end{equation}
be the distribution function of functional (\ref{E:fun}) 
 and let $u(p,q)$ denote its double Laplace--Stieltjes transform:
\begin{equation} \label{E:Laplace}
	u(p,q) = \int_{0}^{+\infty}\int_{0}^{+\infty} e^{-p\lambda - qt}
	d_{\lambda} F(\lambda,t)dt.
\end{equation}

%
 A direct application of the classical result of Kac (1949), see also Rosenblatt (1951), shows
that the function  $u(p,q)$ can be recovered in the form 
\begin{equation} 
	u(p,q) =  \int_{-\infty }^{+\infty} u(p,q,x) \,dx,
\end{equation}
where $u(p,q,x)$ is a solution to the following second order ordinary differential equation
%
%
\begin{equation} \label{E:initial}
	u''(p,q,x) - 2 (q + p e^{2x}) \, u(p,q,x) = 0,  \text{  where $x \ne 0$}, 
\end{equation}
for which we seek a continuous solution subject to conditions
\begin{equation} 
\begin{align}                                      \label{E:initial1}  
	&u(p,q,x) \to 0,  \text{ as $x \to \pm \infty$}, \\  \label{E:initial1a} 
	&|u'(p,q,x)| < M,  \text{  $x \ne 0$}, \\
	&u'(p,q,0+) - u'(p,q,0-) = -2.                       \label{E:initial2}  
\end{align}
\end{equation}
The distribution $F(\lambda,t)$ and its density $f(\lambda,t) :=  F_{\lambda}(\lambda,t)$
can be recovered from $u(p,q)$ by double Laplace inversions:
\begin{equation} \label{E:LaplaceInva}
	F(\lambda,t)=\frac{1}{(2\pi i)^{2}} 
	\int_{a-i \infty}^{a+i \infty}
 	dq e^{qt}\int_{b-i \infty}^{b+i \infty}\frac{1}{p} \, u(p,q) \, 
 	e^{p\lambda}dp
\end{equation}
\begin{equation} \label{E:LaplaceInvb}
	f(\lambda,t)=\frac{1}{(2\pi i)^{2}} 
	\int_{a-i \infty}^{a+i \infty}
	 dq e^{qt}\int_{b-i \infty}^{b+i \infty} u(p,q) \, e^{p\lambda}dp
\end{equation}
with any $ a>0, b>0 $. \\
%
%
%
%
\begin{Thm} \label{Th:1} 
The general solution of the equation
\begin{equation} \label{E:firstDist}
y''(x) - 2 (q + p\, e^{2x}) \, y(x) = 0
\end{equation} 
is given by
\begin{equation} \label{E:firstDist}
y(x) = c_1 I_{\sqrt{2q}} (\sqrt{2p} \, e^{x}) + c_2 K_{\sqrt{2q}} (\sqrt{2p}\, e^{x}), 
\end{equation}
where $I_{\nu}(x)$ and $K_{\nu}(x)$ are modified Bessel functions of the first and second kind.
\end{Thm}
\begin{pf}
Equation (\ref{E:initial}) is a particular case of the equation
\begin{equation} \notag
x^2 y'' - 2a xy' + [(b^2 e^{2cx} - \nu^2) c^2x^2 + a(a-1)]y = 0,
\end{equation} 
where $a=0$, $c=1$, $b^2=-2p$, $\nu^2=2q$, see Kamke [1959] 2.163 (23),
and its solution is
\begin{equation} \notag
\begin{align}    \notag
y &= x^{-a} Z_{\nu}(b e^{cx}) = Z_{\sqrt{2q}}(\sqrt{- 2p} \, e^{x}) = Z_{\sqrt{2q}}(\pm i \sqrt{2p}\, e^{x}) = \\
&= c_1 I_{\sqrt{2q}} (\sqrt{2p} \, e^{x}) + c_2 K_{\sqrt{2q}} (\sqrt{2p}\, e^{x}). \notag
\end{align}  
\end{equation}  
\end{pf}
\begin{Thm}
The solution of equation (\ref{E:initial}) subject to conditions  (\ref{E:initial1}) -- (\ref{E:initial2}) is of the form
\begin{equation} \notag
u(p,q,x) = 
\begin{cases}    \notag
& 2 I_{\sqrt{2q}} (\sqrt{2p})  K_{\sqrt{2q}} (\sqrt{2p}\, e^{x}), \text{ if } x \ge 0, \\
&2 K_{\sqrt{2q}} (\sqrt{2p})  I_{\sqrt{2q}} (\sqrt{2p}\, e^{x}),   \text{ if } x < 0 .
\end{cases}  
\end{equation}
\end{Thm}
\begin{pf}
According to Theorem \ref{Th:1} we seek a solution of the form
\begin{equation} \notag
u(p,q,x) = 
\begin{cases}    \notag
& c_1 I_{\sqrt{2q}} (\sqrt{2p} \, e^{x}) + c_2 K_{\sqrt{2q}} (\sqrt{2p}\, e^{x}), \text{ if } x \ge 0, \\
&c'_1 I_{\sqrt{2q}} (\sqrt{2p} \, e^{x}) + c'_2 K_{\sqrt{2q}} (\sqrt{2p}\, e^{x}),   \text{ if } x < 0,
\end{cases}  
\end{equation}
where constants $c_1, \ldots, c'_2$ to be determined by conditions  (\ref{E:initial1}) -- (\ref{E:initial2}). \\
Since $K_{\nu}(x) \to \infty$ as $x \to 0+$ and it  decays as $x \to \infty$, while $I_{\nu}(x) \to 0 $ as $x \to 0+$ and  $I_{\nu}(x) \to \infty$   as $x \to \infty$, in order to satisfy condition (\ref{E:initial1}) we set  
\[
c_1 = c'_2 = 0,
\]
hence
\begin{equation} \notag
u(p,q,x) = 
\begin{cases}    \notag
& c_2 K_{\sqrt{2q}} (\sqrt{2p}\, e^{x}), \text{ if } x \ge 0, \\
&c'_1 I_{\sqrt{2q}} (\sqrt{2p} \, e^{x}),   \text{ if } x < 0,
\end{cases}  
\end{equation}
and this function obviously satisfies (\ref{E:initial1a}). \\

It remains to show that the choice
\begin{equation}                         \notag
\begin{align}                            \notag
&c_1 = 2 I_{\sqrt{2q}} (\sqrt{2p}) , \\ \notag
&c'_2 = 2 K_{\sqrt{2q}} (\sqrt{2p})
\end{align}
\end{equation} 
generates a continuous solution that satisfies (\ref{E:initial2}). \\
The continuity follows from
\begin{equation}                         \notag
u(p,q,0+) = u(p,q,0-) = 
2 \,I_{\sqrt{2q}} (\sqrt{2p})  K_{\sqrt{2q}} (\sqrt{2p}). 
\end{equation} 

By a straightforward computation we obtain: 
\begin{equation}                         \notag
\begin{align}                            \notag
&u'(p,q,0+) - u'(p,q,0-) = \\
&2 \sqrt{2p}\,\big[I_{\sqrt{2q}} (\sqrt{2p})  K'_{\sqrt{2q}} (\sqrt{2p}) - I'_{\sqrt{2q}} (\sqrt{2p})  K_{\sqrt{2q}} (\sqrt{2p}) \big] 
\end{align}
\end{equation} 
We observe that the expression in the brackets is up to the sign the Wronskian: 
\[
-W\{K_{\nu}(z), I_{\nu}(z) \} = -\frac{1}{z},
\]
see Olver (1974) p.251 (8.07), where $z= \sqrt{2p} $,  $\nu = \sqrt{2q} $, and this yields
\[
u'(p,q,0+) - u'(p,q,0-) = -2.
\]
\end{pf}
In order to recover the density distribution $f(\lambda,t)$ from $u(p,q,x)$ we first invert the latter with respect to $p$ (Section \ref{S:pinv}), then integrate the result on $x$ (Section \ref{S:intx} ), and then invert with respect to $q$ in the following sections. 
These steps involve interchange of integrations which is easily justifiable by simple estimates and Fubini's theorem,
see for example Tolmatz (2000), (2003). 
%
\section{Laplace Inversion with Respect to $p$} \label{S:pinv}
The following theorem provides a Laplace inversion of $u(p,q,x)$.
\begin{Thm}
\begin{equation} \label{E:guess}
w(\lambda,q,x) := \frac{1}{2\pi i} \int_{c-i\infty}^{c+i\infty} u(p,q,x) \, e^{p \lambda} \,dp = 
\frac{1}{\lambda} \exp \bigg(-\frac{e^{2x}+1}{2\lambda} \bigg) \, I_{\sqrt{2q}} \bigg(\frac{e^x}{\lambda} \bigg). 
\end{equation}
\end{Thm}
\begin{pf}
It suffices to verify that the Laplace transform of the right--hand side of (\ref{E:guess}) results in $u(p,q,x)$, that is  
\begin{equation} \label{E:firstDist}
\int_0^\infty w(\lambda,q,x) e^{-p\lambda} \,d\lambda = \int_{0}^{\infty} 
\frac{1}{\lambda} \exp \bigg(-\frac{e^{2x}+1}{2\lambda} \bigg) \, I_{\sqrt{2q}} \bigg(\frac{e^x}{\lambda} \bigg) e^{-p \lambda} \, d\lambda = u(p,q,x). 
\end{equation}
To this end we make use of a table of Laplace transforms, see for example 
Ditkin and Prudnikov (1961) $9.351$: 
\begin{equation} \label{E:table}
 \int_{0}^{\infty} 
\frac{1}{t} \exp \bigg(-\frac{\alpha + \beta}{2t} \bigg) \, I_{\nu} \bigg(\frac{\alpha - \beta}{2t} \bigg) e^{-p t} \, dt = 
2K_{\nu}\big[(\sqrt\alpha + \sqrt\beta)\sqrt p \, \big]\, I_{\nu}\big[(\sqrt\alpha - \sqrt\beta)\sqrt p \, \big], 
\end{equation}

where $\Re \alpha \ge \Re \beta >0.$ \\
We apply (\ref{E:table}) with 
\begin{equation} \label{E:a1}    \notag                  
\begin{cases}
&\alpha + \beta = e^{2x}+1, \\ 
&\alpha - \beta = 2e^x,
\end{cases} 
\end{equation} 
that is
\begin{equation} \label{E:a2}    \notag 
\begin{cases}                           
&\alpha = \frac{1}{2} (e^x + 1)^2, \\ 
& \beta =\frac{1}{2} (e^x - 1)^2 , 
\end{cases}
\end{equation} 
\begin{equation}  
\begin{cases}                            \notag
&\sqrt\alpha = \frac{1}{\sqrt 2} \, (e^x + 1), \\ 
& \sqrt\beta =\frac{1}{\sqrt 2} \, |e^x - 1| , 
\end{cases}
\end{equation} 
and obtain
\begin{equation}      \label{E:a3}           
\sqrt\alpha + \sqrt\beta =  
\begin{cases}                           
&\sqrt 2 \, e^x , \text{ if } x \ge 0, \\ 
& \sqrt 2, \text{ if } x  < 0, 
\end{cases}
\end{equation} 
\begin{equation}                 \label{E:a4}
\sqrt\alpha - \sqrt\beta =  
\begin{cases}                        
&\sqrt 2, \text{ if } x \ge 0, \\ 
& \sqrt 2 \, e^x , \text{ if } x  < 0. 
\end{cases}
\end{equation} 
Substitution of (\ref{E:a3})--(\ref{E:a4}) in (\ref{E:table}) yields    
\begin{equation} \label{E:firstDist}
\begin{align}   \notag
 \int_{0}^{\infty} 
&\frac{1}{\lambda} \exp \bigg(-\frac{e^{2x}+1}{2\lambda} \bigg) \, I_{\sqrt{2q}} \bigg(\frac{e^x}{\lambda} \bigg) e^{-p \lambda} \, d\lambda = \\ \notag 
&=
\begin{cases}                        
&2 I_{\sqrt{2q}} (\sqrt{2p})  K_{\sqrt{2q}} (\sqrt{2p}\, e^{x}), \text{ if } x \ge 0, \\ 
&2 K_{\sqrt{2q}} (\sqrt{2p})  I_{\sqrt{2q}} (\sqrt{2p}\, e^{x}), \text{ if } x  < 0. 
\end{cases}
\\ \notag 
&= u(p,q,x).
\end{align} 
\end{equation}
\end{pf}
{\em Remark}. Alternatively, we can prove the above theorem by a straightforward derivation without "guessing" the function in (\ref{E:guess}), however that proof involves tedious computations. \\
{\em Remark}. $w(\lambda,q,x)$ is a Laplace transform of the joint density distribution of random variables $\int_0^t e^{2 B(s)} ds$ and $B(t)$.  

\section{The Integral $\int w\,dx$} \label{S:intx}
%
In this section we compute the integral 
\begin{equation} \label{E:firstDist}
w(\lambda,q) := \int_{-\infty}^{\infty} u(\lambda,q,x) \, dx.  \notag
\end{equation}
The following proposition will be used in the proof of the next theorem.
\begin{Prop} \label{P:GR1}
Let $\Re \alpha > 0$, $\Re {(\mu + \nu)} > -1$, $\beta >0 $. Then
\begin{equation} 
\begin{align}
  & \int_{0}^{\infty} x^{\mu} e^{-\alpha x^2} I_{\nu}(\beta x) \, dx = \\  \notag
& \frac{\beta^\nu \Gamma \big( \frac{1}{2}\nu +  \frac{1}{2}\mu + \frac{1}{2} \big)}
{2^{\nu+1} \alpha^{\frac{1}{2}(\mu+\nu+1)} \Gamma (\nu+1)}\,
{}_1F_{1} \bigg( \frac{\nu+\mu+1}{2}; \nu+1; \frac{\beta^2}{4\alpha}    \bigg)
\end{align}
\end{equation}
\end{Prop}
\begin{pf}
By Gradshtein and Ryzhik (1963) {\bf 6.631} 1.
\begin{equation}  \label{E:anal}  
\begin{align}     
 	& \int_{0}^{\infty} x^{\mu} e^{-\alpha x^2} J_{\nu}(\beta x)\, dx = \\ \notag 
	&\frac{\beta^\nu \Gamma \big( \frac{1}{2}\nu +  \frac{1}{2}\mu + \frac{1}{2} \big)}
	{2^{\nu+1} \alpha^{\frac{1}{2}(\mu+\nu+1)} \Gamma (\nu+1)}\,
	{}_1F_{1} \bigg( \frac{\nu+\mu+1}{2}; \nu+1; -\frac{\beta^2}{4\alpha}    \bigg)
\end{align}
\end{equation}
The integral in (\ref{E:anal}) is a function of the real parameter $\beta$ that admits an analytic continuation into the domain $D:=\{\beta | \arg \beta < \pi  \}$ in the complex $\beta$--plane. \\
Indeed, for any fixed $x>0$ 
$J_{\nu}(\beta x)$ is analytic in $D$, while the asymptotics of $J_{\nu}(z)$ in the domain $|\arg z|<\pi$ according to Gradshtein and Ryzhik (1963) {\bf 8.451} 1. is

\begin{equation} \label{E:anal1} \notag
 J_{\nu}(z) = \sqrt{\frac{2}{\pi z}}\bigg[\cos{\bigg(z - \frac{\pi}{2}\nu - \frac{\pi}{4}
   \bigg)} \big(1 + O( 1/|z|^2) \big) 
\bigg]   \text{ as } |z| \to \infty,
\end{equation}
%
and this implies
\begin{equation}  \notag
 J_{\nu}(\beta x) = \sqrt{\frac{2}{\pi z}}\bigg[\cos{\bigg(\beta x - \frac{\pi}{2}\nu - \frac{\pi}{4}
   \bigg)} \big(1 + O( 1/|\beta|^2) \big) 
\bigg]   \text{ as } |\beta| \to \infty,
\end{equation}
which ensures that for any $\beta \in D$ the integral in (\ref{E:anal}) converges and is differentiable on $\beta$.\\ 
The function in the right--hand side of (\ref{E:anal}) is also analytic in $D$ due to the analycity of
 ${}_1F_{1}(\alpha; \gamma; z) $   in the whole $z$--plane 
and since $\beta^\nu$ is analytic in $D$.
Therefore (\ref{E:anal}) holds also on the imaginary $\beta$-axis, and we get 
\begin{equation} \label{E:anal1} \notag
 \int_{0}^{\infty} x^{\mu} e^{-\alpha x^2} J_{\nu}(i\beta x)\, dx = 
\frac{(i\beta)^\nu \Gamma \big( \frac{1}{2}\nu +  \frac{1}{2}\mu + \frac{1}{2} \big)}
{2^{\nu+1} \alpha^{\frac{1}{2}(\mu+\nu+1)} \Gamma (\nu+1)}\,
{}_1F_{1} \bigg( \frac{\nu+\mu+1}{2}; \nu+1; -\frac{(i\beta)^2}{4\alpha}    \bigg)
\end{equation}
The proof is completed by substituting
\begin{equation} \notag
 J_{\nu}(ix) = i^{\nu}I_{\nu}(x)
\end{equation}
in the previous relationship and obvious simplifications. 

\end{pf}
\begin{Thm}\label{Th:w}
\begin{equation} 
w(\lambda,q) =\frac{1}{2\lambda} \,  e^{-1/(2\lambda)}(2\lambda)^{-\sqrt {2q}/2} \frac{\Gamma\big( \frac{\sqrt {2q}}{2} \big)}{\Gamma(\sqrt {2q}+1)} \, { }_{1}F_{1} \bigg( \frac{\sqrt {2q}}{2}\, ; \sqrt {2q}+1; \frac{1}{2\lambda}  \bigg). 
%
\end{equation}
\end{Thm}
\begin{pf}
In the integral
\begin{equation} \notag
w(\lambda,q) = \int_{-\infty}^{\infty} u(\lambda,q,x) \, dx =\int_{-\infty}^{\infty} \frac{1}{\lambda} \exp \bigg(-\frac{e^{2x}+1}{2\lambda} \bigg) \, I_{\sqrt{2q}} \bigg(\frac{e^x}{\lambda} \bigg) \, dx 
\end{equation}
we perform the change of variables
\begin{equation} \notag
\begin{align} \notag
&y = e^x, \\ \notag
&x = \ln y, \\ \notag
&dx = \frac{dy}{y},
\end{align}
\end{equation}
and obtain by Proposition \ref{P:GR1} with $\mu= -1$, $\alpha= 1/(2\lambda)$, $\beta=1/\lambda$, $\nu = \sqrt{2q}$

\begin{equation} \label{E:firstDist}
\begin{align}
	&w(\lambda,q) = \frac{1}{\lambda} \,  e^{-1/(2\lambda)}\int_{-\infty}^{\infty}  \exp \bigg(-\frac{e^{2x}}{2\lambda} \bigg) \, I_{\sqrt{2q}} \bigg(\frac{e^x}{\lambda} \bigg) \, dx =  \\ \notag
	&\frac{1}{\lambda} \,  e^{-1/(2\lambda)}\int_0^{\infty} 
	\frac{1}{y} \, e^ {-y^2/(2\lambda)} \, I_{\sqrt{2q}} \bigg(\frac{y}{\lambda} \bigg) \, dy = \\ \notag
	&\frac{1}{2\lambda} \,  e^{-1/(2\lambda)}(2\lambda)^{-\sqrt {2q}/2} \frac{\Gamma\big( \frac{\sqrt {2q}}{2} \big)}{\Gamma(\sqrt {2q}+1)} \, { }_{1}F_{1} \bigg( \frac{\sqrt {2q}}{2}\, ; \sqrt {2q}+1; \frac{1}{2\lambda}  \bigg). 
\end{align}
\end{equation}
\end{pf}

In order to compute the distribution density $f(\lambda,t)$ of the functional $\int_0^t e^{2 B(s)} ds$ we seek to invert $w(\lambda, q)$, that is to compute the integral
\begin{equation} 
\begin{align}
f(\lambda,t) &= \frac{1}{2\pi i} \int_{c-i\infty}^{c+i\infty} w(\lambda, q)  e^{q t} \, dq = \\  \notag
&= \frac{1}{2 \lambda} \,  e^{-1/(2\lambda)} \frac{1}{2\pi i}\int_{c-i\infty}^{c+i\infty}(2\lambda)^{-\sqrt {2q}/2} \frac{\Gamma\big( \frac{\sqrt {2q}}{2} \big)}{\Gamma(\sqrt {2q}+1)} \, { }_{1}F_{1} \bigg( \frac{\sqrt {2q}}{2}\, ; \sqrt {2q}+1; \frac{1}{2\lambda}  \bigg) 
\, e^{qt} \,dq. 
\end{align}
\end{equation}
These computations will be performed in Sections \ref{S:w}--\ref{S:lt}.
%

\section{A Series for $w(\lambda, q)$} \label{S:w}
%
\begin{Thm} \label{Th:w1}
\begin{equation} \label{E:wseries}
w(\lambda, q) = \frac{1}{2\lambda} \,  e^{-1/(2\lambda)} \sum_{n=0}^{\infty} \frac{1}{n!}\,
\, w_n(\lambda, q)\frac{1}{(2\lambda)^n},
\end{equation}
where  
\begin{equation}
w_n(\lambda, q) := (2\lambda)^{-\sqrt {2q}/2} \frac{\Gamma\big( \frac{\sqrt {2q}}{2}+n \big)}{\Gamma(\sqrt {2q}+n+1)}. 
\end{equation}
\end{Thm}
\begin{pf}
The proof follows from Theorem \ref{Th:w}, the expansion
\begin{equation}  \notag
{}_1F_1(\alpha; \gamma;z) = 1 + \frac{\alpha}{\gamma}\, \frac{z}{1!} + \frac{\alpha (\alpha+1)}{\gamma (\gamma+1)}\, \frac{z^2}{2!} \ldots ,
\end{equation}
where  $\alpha = \sqrt{2q}/2 $, $\gamma = \sqrt{2q}+1$, $z = 1/(2\lambda)$ and the identity
\begin{equation} \notag
\Gamma(a+n) = \Gamma(a)\, a \ldots (a+n-1),
\end{equation}
where $a = \sqrt{2q}/2, \, \sqrt{2q}+1$.
\end{pf}
In the remaining part of this section we shall justify a term--wise inversion of the series in  (\ref{E:wseries}). \\
The following elementary lemma is given without proof.
\begin{Lem}
Let $|\arg z| \le  \pi/2$. Then for any  $s \ge 0$  
\begin{equation}
\bigg|\frac{z/2 + s}{z + s} \bigg| < 1.
\end{equation}
\end{Lem}
\begin{Cor}\label{C:est1}
For $ |\arg q |< \pi $, $n = 0,1, \ldots $ holds
\begin{equation}
\bigg|\frac{\sqrt{2q}/2 + n}{\sqrt{2q} + n} \bigg| < 1.
\end{equation}
\end{Cor}
\begin{Lem}
Let $|\Re q| > 0 $,  $n = 0,1, \ldots $  . \\
Then  
\begin{equation}
\bigg|\frac{\Gamma(\sqrt{2q}/2 + n)}{\Gamma(\sqrt{2q} + n +1)} \bigg| < \frac{1}{|\sqrt{2q}|}\bigg|\frac{\Gamma(\sqrt{2q}/2 )}{\Gamma(\sqrt{2q})} \bigg|.
\end{equation}
\end{Lem}
\begin{pf}
The proof follows by elementary transforms and Corollary \ref{C:est1}:
\begin{equation}
\begin{align} \notag
&\bigg|\frac{\Gamma(\sqrt{2q}/2 + n)}{\Gamma(\sqrt{2q} + n +1)} \bigg| =
\bigg|\frac{1}{\sqrt{2q}+n} \cdot
\frac{(\sqrt{2q}/2+n-1) \ldots (\sqrt{2q}/2)\,\Gamma(\sqrt{2q}/2 )}
{(\sqrt{2q}+n-1) \ldots (\sqrt{2q})\,\Gamma(\sqrt{2q})} \bigg| < \\ \notag
& <\bigg|\frac{1}{\sqrt{2q}+n} 
 \frac{\Gamma(\sqrt{2q}/2 )} {\Gamma(\sqrt{2q})} \bigg| <  \frac{1}{|\sqrt{2q}|}\bigg|\frac{\Gamma(\sqrt{2q}/2 )}{\Gamma(\sqrt{2q})} \bigg|. 
\end{align}
\end{equation}
\end{pf}
\begin{Lem}
Let $|\Re q| > 0 $,  $n = 0,1, \ldots $  . \\
Then  
\begin{equation}  \label{E:gest}
\bigg|\frac{\Gamma(\sqrt{2q}/2 + n)}{\Gamma(\sqrt{2q} + n +1)} \bigg| < 2\sqrt \pi \,\frac{1}{|\sqrt{2q}|}\big|2^{-\sqrt{2q}} \,\big| \,
\bigg|\frac{1}{\Gamma(\sqrt{2q}/2 +1/2)} \bigg|.
\end{equation}
\end{Lem}

\begin{pf}
The proof follows by the previous lemma and Legendre's formula   
\begin{equation}
\Gamma(2a) = \frac{1}{\sqrt \pi} 2^{2a-1} \Gamma(a) \Gamma(a + 1/2),  
\end{equation}
which we apply in the form
\begin{equation}
\frac{\Gamma(a)}{\Gamma(2a)} = 2\sqrt\pi \, 2^{-2a}\frac{1}{\Gamma(a + 1/2)},  
\end{equation}
where $a=\sqrt{2q}/2$.
\end{pf}
\begin{Lem}
Let $q(\tau) = c +i\tau$, where $c>0$ and $-\infty < \tau < \infty$ is a parameter. \\
Then
\begin{equation}
\bigg|\frac{\Gamma(\sqrt{2q(\tau)}/2 + n)}{\Gamma(\sqrt{2q(\tau)} + n +1)} \bigg| < C \, |\tau|^{-a \sqrt{|\tau|}},
\end{equation}
where $C, a >0$ are some constants not depending on $n$.
\end{Lem}

\begin{pf}
The proof readily follows from the previous lemma and Stirling's formula 
\begin{equation} \label{E:stirling}
\Gamma(z) \thicksim e^{-z} z^z \bigg( \frac{2\pi}{z} \bigg)^{1/2}
\end{equation}
as $z \to \infty$ in the sector $|\arg z|\le \pi-\delta $, which implies
\begin{equation}
\frac{1}{\Gamma \big(\sqrt{2(c+i\tau)}/2 +1/2 \big)} =  O \bigg(|\tau|^{-a \sqrt{|\tau|}} \bigg)
\end{equation}
as $\tau \to \pm \infty$.

\end{pf}
The following theorem shows that the series in Theorem \ref{Th:w1} can be inverted term by term.
\begin{Thm}\label{Th:fser}
\begin{equation}
f(\lambda, t) = \frac{1}{2\lambda} \,  e^{-1/(2\lambda)} \sum_{n=0}^{\infty} \frac{1}{n!}\,
\, f_n(\lambda, t)\frac{1}{(2\lambda)^n},
\end{equation}
where  
\begin{equation}\label{E:fn0}
 f_n(\lambda, t) := \frac{1}{2\pi i} \int_{c-i\infty}^{c + i\infty} w_n(\lambda, q) \, e^{qt}\,dq,    
\end{equation}
where
\begin{equation}\label{E:wn0} 
w_n(\lambda, q) := (2\lambda)^{-\sqrt {2q}/2} \frac{\Gamma\big( \frac{\sqrt {2q}}{2}+n \big)}{\Gamma(\sqrt {2q}+n+1)}, 
\end{equation}
and this series converges absolutely for any fixed $\lambda, t > 0$. 
\end{Thm}
\begin{pf}
The previous lemma implies absolute summability of 
\begin{equation}
\frac{1}{n!}(2\lambda)^{-\sqrt {2q(\tau)}/2} \frac{\Gamma\big( \frac{\sqrt {2q(\tau)}}{2}+n \big)}{\Gamma(\sqrt {2q(\tau)}+n+1)} \, e^{q(\tau)t}\frac{1}{(2\lambda)^n} 
\end{equation}
as a function of variables $\tau$, $n$ with any fixed $\lambda, t > 0$ in the sense 
\begin{equation}
\sum_{n=0}^{\infty} \frac{1}{n!}\int_{-\infty}^{\infty} \bigg|(2\lambda)^{-\sqrt {2q(\tau)}/2} \frac{\Gamma\big( \frac{\sqrt {2q(\tau)}}{2}+n \big)}{\Gamma(\sqrt {2q(\tau)}+n+1)} e^{q(\tau)t} \bigg| \frac{1}{(2\lambda)^n} d\tau < \infty,
\end{equation}

and the proof follows by Fubini's theorem.
\end{pf}

\section{Transforms of $w_n(\lambda, q)$ } \label{S:ft}
%
We wish to prepare the integrals
\begin{equation}
 f_n(\lambda, t) = \frac{1}{2\pi i} \int_{c-i\infty}^{c+i\infty} w_n(\lambda, q) \, e^{qt}\,dq    
\end{equation}
for an explicit computation.
\begin{Prop}\label{P:w}
\begin{equation} \label{E:pw}
\begin{align}
 	 f_n(\lambda, t) &= \frac{1}{2\pi i} \int_{c-i\infty}^{c+i\infty} w_n(\lambda, q) \, e^{qt}\,dq  = \\ \notag
	&=\frac{e^{t/2}}{\sqrt{2 \lambda}}\frac{1}{2\pi i} \int_{c-i\infty}^{c+i\infty}
(2\lambda e^{-2t})^{-\sqrt z} \frac{\Gamma(\sqrt z+n+1/2)}{\Gamma(2\sqrt z+n+2)}\, 
	\frac{2\sqrt z + 1}{\sqrt z}
 	e^{2zt}\,dz.
\end{align} 
\end{equation} 
\end{Prop}

\begin{pf}
The proof follows by the reparametrization $\sqrt{2q} = 2 \sqrt{z}+1$, $q=2z + 2 \sqrt{z}+1/2 $,
$dq = (2\sqrt z + 1)/(\sqrt z)\,dz$.
\end{pf}

We shall need auxiliary functions 
\begin{equation}
\varphi_n(w):= 2^{2w}\frac{\Gamma(w+n+1/2)}{\Gamma(2w+n+2)} 
\frac{2w + 1}{w}, 
\end{equation} 
where $n = 0, 1, \ldots$. \\
By making use of the properties of the Gamma function as in the previous section we can write
\begin{equation}
\varphi_0(w)= \frac{2^{2w}}{w\Gamma(w+1)}, 
\end{equation} 
while for $n=1,2, \ldots$
\begin{equation} \label{E:ratnl}
\begin{align}
&\varphi_n(w)=2^{2w} \frac{\Gamma(w+n+1/2)}{\Gamma(2w+n+2)} 
\frac{2w + 1}{w} = \\ \notag
& \sqrt\pi \,  2^{-(n+1)} \frac{1}{w+n/2+1/2}\,\frac{(w+n-1/2)\cdot \ldots \cdot (w+1/2)}
{(w+n/2)\cdot \ldots \cdot (w+1/2)} \,\frac{2w + 1}{w} \frac{1}{\Gamma(w+1)} = \\ \notag
& \sqrt\pi \,  2^{-(n+1)} \frac{1}{w+n/2+1/2}\, \bigg( \prod_{j=1}^n \frac{w+j-1/2}{w+j/2} \bigg)  \frac{2w + 1}{w}\frac{1}{\Gamma(w+1)} .   
\end{align}
\end{equation} 
Consider the rational factor in (\ref{E:ratnl})

\begin{equation}
r_n(w):=
\sqrt\pi \,  2^{-(n+1)} \frac{1}{w+n/2+1/2}\, \bigg( \prod_{j=1}^n \frac{w+j-1/2}{w+j/2} \bigg)  \frac{2w + 1}{w}.
\end{equation} 
Observe that due to cancellations of certain factors in the nominator and denominator $r_0(w)$ is regular except the simple pole at $w=0$, and $r_n(w)$, $n \ge 1$ is regular except the simple poles at $w = 0, -1, \ldots, -\big[\frac{n-1}{2} \big] - 1$.

We wish to decompose $r_n(w)$ into a sum of elementary fractions.
\begin{Thm}\label{Th:ank}
\begin{equation} \label{E:decomp}
r_n(w)= \sum_{k=0}^{[(n-1)/2]+1} \frac{a_k^{(n)}}{w+k},
\end{equation}
where
\begin{equation}
a_0^{(n)}  = \frac{\Gamma(n+1/2)}{\Gamma(n+2)}, 
\end{equation} 

\begin{equation}
a_k^{(n)}  = 2^{-2k}\frac{(-1)^{k-1}}{(k-1)!} \frac{\Gamma(n -k +1/2)}{\Gamma(n - 2k +2)}, 
\end{equation}
where $k = 1, \ldots, [(n-1)/2]+1$. 
\end{Thm}
\begin{pf}
Since all the poles of $r_n(w)$ are simple, the existence of decomposition (\ref{E:decomp}) is a well--known fact. \\
%
The relationship 
\begin{equation} \label{E:rat}
\begin{align}
&\varphi_n(w)\Gamma(w+1)=2^{2w} \frac{\Gamma(w+n+1/2)}{\Gamma(2w+n+2)} 
\frac{2w + 1}{w} = \\ \notag
& \sqrt\pi \,  2^{-(n+1)} \frac{1}{w+n/2+1/2}\,\frac{(w+n-1/2)\cdot \ldots \cdot (w+1/2)}
{(w+n/2)\cdot \ldots \cdot (w+1/2)} \,\frac{2w + 1}{w}= \sum_{k=0}^{[(n-1)/2]+1} \frac{a_k^{(n)}}{w+k}
\end{align}
\end{equation} 
shows that the coefficients $a_k^{(n)}$ are the residues of $\varphi_n(w)\Gamma(w+1)$ at $w=0, -1, \ldots, -[(n-1)/2]-1$, hence
\begin{equation}
\begin{align}    \notag
a_k^{(n)}  &= \text{\em Res}_{w=-k}\{\varphi_n(w)\Gamma(w+1) \}=
\text{\em Res}_{w=-k}  \bigg\{2^{2w} \frac{\Gamma(w+n+1/2)}{\Gamma(2w+n+2)} 
\frac{2w + 1}{w} \Gamma(w+1)\bigg\}= \\
&=
\begin{cases}
&\frac{\Gamma(n+1/2)}{\Gamma(n+2)}, \text{ if } k=0,  \\  \notag
&\varphi_n(-k) \text{\em Res}_{w=-k} \Gamma(w+1), \text{ if } k=1,2, \ldots, [\frac{n-1}{2} ]+1,
\end{cases}
 \\
&=
\begin{cases}
&\frac{\Gamma(n+1/2)}{\Gamma(n+2)}, \text{ if } k=0,  \\  \notag
&\varphi_n(-k) \frac{(-1)^{k-1}}{(k-1)!}, \text{ if } k=1,2, \ldots, [\frac{n-1}{2} ]+1,     
\end{cases}
 \\
&=
\begin{cases}
&\frac{\Gamma(n+1/2)}{\Gamma(n+2)}, \text{ if } k=0,  \\  \notag
&2^{-2k}\frac{(-1)^{k-1}}{(k-1)!} \frac{\Gamma(n -k +1/2)}{\Gamma(n - 2k +2)}, \text{ if } k=1,2, \ldots, [\frac{n-1}{2} ]+1.
\end{cases}
\end{align}
\end{equation} 
\end{pf}
%
For typographical reasons we introduce the notation
\begin{equation} \label{E:fseries}
a := 8\lambda e^{-2t}.
\end{equation}
\begin{Thm}\label{Th:laplace2}
\begin{equation} \label{E:fseries}
	f(\lambda, t) = \frac{1}{2\lambda \sqrt{2\lambda}} \,  e^{-1/(2\lambda)} e^{t/2}\sum_{n=0}^{\infty} \frac{1}{n!}\,
	\, f_n(\lambda, t)\frac{1}{(2\lambda)^n},
\end{equation}

where
\begin{equation}
 	f_0(\lambda, t) = \frac{\sqrt \pi}{2\pi i} \int_{c-i\infty}^{c+i\infty} \frac{a^{-\sqrt z} \, e^{2zt}}{\sqrt z \, \Gamma(\sqrt z +1)} \, dz,    
\end{equation}
and for $n = 1,2, \ldots$
\begin{equation}
 f_n(\lambda, t) = \sum_{k=0}^{[(n-1)/2]+1} a_k^{(n)}  \frac{1}{2\pi i} 
\int_{c-i\infty}^{c+i\infty}\frac{a^{-\sqrt z} \, e^{2zt}}{(\sqrt z+k) \, \Gamma(\sqrt z +1)} \, dz,     
\end{equation}
where coefficients $a_k^{(n)}$ are defined in Theorem \ref{Th:ank}.
\end{Thm}
\begin{pf}
Only the case $n \ge 1$ needs a proof. By Theorem \ref{Th:ank} we can write
\begin{equation}  \notag
	\begin{align} \notag
	\frac{\Gamma(\sqrt z+n+1/2)}{\Gamma(2\sqrt z+n+2)}\, \frac{2\sqrt z + 1}{\sqrt z} &=  2^{-2\sqrt z} \varphi_n(\sqrt z) =  
	2^{-2\sqrt z} \frac{r_n(\sqrt z)}{\Gamma(\sqrt z + 1)} = \\ \notag 
	&=2^{-2\sqrt z} \frac{1}{\Gamma(\sqrt z + 1)} \sum_{k=0}^{[(n-1)/2]+1} \frac{a_k^{(n)}}{\sqrt z + k} ,
	\end{align}
\end{equation}
substituting of which in (\ref{E:pw}) yields the proof.
\end{pf} 
%
\section{Auxiliary Results} \label{S:aux}
%
The following expansion is given in Gradshtein and Ryzhik (1963) {\bf 8.321} 2:
\begin{equation}\label{E:dm1}
	\frac{1}{ \Gamma(z +1)} = \sum_{k=0}^{\infty} d_k \,  z^k,   
\end{equation}
where
\begin{equation}\label{E:dm2} 
\begin{cases}
	&d_0 = 1, \\
	&d_{n+1} = \frac{1}{n+1} \sum_{k=0}^n (-1)^k s_{k+1} d_{n-k},   
\end{cases}
\end{equation}
where $s_1 = {\bf C}$, $s_n = \zeta(n)$ for $n \ge 2$. 
\begin{Thm}
\begin{equation} \label{E:d}
	|d_n| = O(n^{-(1-\varepsilon) n} ).
\end{equation}
\end{Thm}
\begin{pf}
By Cauchy's formula for Taylor coefficients we can write
\begin{equation}  \notag
	d_n = \frac{1}{2\pi i}\int_{|z|=1}  \frac{1}{z^{n+1} \Gamma(z +1)} \,dz, 
\end{equation}
and then estimate this integral by integrating over the circle $z= r_n e^{i\varphi}$, where 
$r_n = n/\ln n$ :   
\begin{equation} \notag
\begin{align}  \notag
	|d_n| &= \bigg| \frac{1}{2\pi i}\int_{|z|=r_n}  \frac{1}{z^{n+1} \Gamma(z +1)} \,dz \bigg| =
 	\frac{1}{2\pi}  \bigg|\int_{|z|=r_n}  \frac{1}{z^{n+2} \Gamma(z)} \,dz \bigg| = \\ \notag
	&= \frac{1}{2\pi}  \bigg|\int_{0}^{2\pi}  \frac{d\varphi}{r_n^{n+1} \, e^{i (n+1) \varphi}\, \Gamma(r_n e^{i\varphi})} \bigg| <   
	\frac{1}{2\pi r_n^{n+1}}  \int_{0}^{2\pi}  \frac{d\varphi} {|\Gamma(r_n e^{i\varphi})|} = \\ \notag
	&=\frac{1}{2\pi r_n^{n+1}}  \int_{-\pi/2}^{\pi/2}  \frac{d\varphi} {|\Gamma(r_n e^{i\varphi})|} +  
	\frac{1}{2\pi r_n^{n+1}}  \int_{\pi /2}^{3/2 \pi}  \frac{d\varphi} {|\Gamma(r_n e^{i\varphi})|}. 
\end{align}
\end{equation}
Estimates for the latter two integrals can be obtained by estimating their integrands by Stirling's formula (\ref{E:stirling}) with making use of Euler's formula
\begin{equation} \label{E:Euler}
	z\Gamma(z) \Gamma(-z) = \frac{\pi}{\sin {\pi z}}
\end{equation}  
in the second integral:
\begin{equation}
\begin{align}  \notag
	&\int_{-\pi/2}^{\pi/2}  \frac{d\varphi} {|\Gamma(r_n e^{i\varphi})|} = \\ \notag 
	&=\int_{-\pi/2}^{\pi/2}  O \bigg(\big| \sqrt {r_n \,e^{i\varphi} } \, e^{r_n\,e^{i\varphi}} (r_n e^{i\varphi})^{-r_n e^{i\varphi}}  \big| \bigg) d\varphi 
	= O \bigg(\sqrt {r_n} \, e^{(1 + \pi/2)  r_n}   \bigg) = 
	 O \bigg( e^{(1 + \pi/2)  r_n}   \bigg), 
\end{align}
\end{equation}
\begin{equation}
\begin{align}  \notag
	&\int_{\pi/2}^{3/2\pi}  \frac{d\varphi} {|\Gamma(r_n e^{i\varphi})|} = \\ \notag
	&=\int_{-\pi/2}^{\pi/2}  O \bigg(\big| \sqrt {r_n e^{i\varphi}}\,\sin{(\pi r_n e^{i\varphi}) } 
	\,e^{r_n\,e^{i\varphi}} (r_n e^{i\varphi})^{r_n e^{i\varphi}}  \big| \bigg) d\varphi 
	= O \bigg(\sqrt{r_n} \,e^{(1 + \pi)  r_n} \,  r_n^{r_n}   \bigg), 
\end{align}
\end{equation}
and therefore
\begin{equation}
\begin{align}  \notag
	|d_n| &= \frac{1}{r_n^{(n+1)}} \bigg\{O \bigg( e^{(1 +\pi/2) \, r_n}   \bigg) +
	O \bigg(\sqrt{r_n} \,e^{(1 + \pi)  r_n} \,r_n^{r_n}   \bigg) \bigg\} =\frac{1}{r_n^{(n+1)}} O \bigg(\sqrt{r_n} \,e^{(1 + \pi)  r_n} \,r_n^{r_n}   \bigg) = \\ \notag 
	&=  O\bigg( e^{(1+\pi) n/ \ln n}  \bigg( \frac{n}{\ln n}   \bigg)^{-n \big(1 - 1/\ln n - 1/(2n) \big)}    \bigg) 
	= O(n^{-(1-\varepsilon) n} ). 
\end{align}
\end{equation}
\end{pf}
{\em Remark}. By applying the saddle point method we can obtain an exact asymptotics for $d_n$ and thus improve the above estimates.
%
\section{Inverse Laplace Transforms} \label{S:lt}
\begin{Thm}
\begin{equation} \label{E:fseries}
	f(\lambda, t) = \frac{1}{2\lambda \sqrt{2\lambda}} \,  e^{-1/(2\lambda)} e^{t/2}\sum_{n=0}^{\infty} \frac{1}{n!}\,
	\, f_n(\lambda, t)\frac{1}{(2\lambda)^n},
\end{equation}
where
\begin{equation}\label{E:f0x} 
 	f_0(\lambda, t) = \frac{\sqrt \pi}{2\pi i} \int_{-\infty}^{(0+)} \frac{a^{-\sqrt z} \, e^{2zt}}{\sqrt z \, \Gamma(\sqrt z +1)} \, dz,    
\end{equation}
and for $n = 1,2, \ldots$
\begin{equation}\label{E:fnx}  
 	f_n(\lambda, t) = \sum_{k=0}^{[(n-1)/2]+1} a_k^{(n)}  \frac{1}{2\pi i} 
	\int_{-\infty}^{(0+)}\frac{a^{-\sqrt z} \, e^{2zt}}{(\sqrt z+k) \, \Gamma(\sqrt z +1)} \, dz,     
\end{equation}
where $a = 8\lambda e^{-2t}$,  and the  coefficients $a_k^{(n)}$ are defined in Theorem \ref{Th:ank}.
\end{Thm}
\begin{pf}
The proof follows from Theorem \ref{Th:laplace2}  by pushing the integration contours to the left, which is justified by Cauchy's theorem and 
obvious estimates by Stirling's formula (\ref{E:stirling}) with making use of Euler's formula (\ref{E:Euler}).
\end{pf} 
\begin{Lem}
\begin{equation}\label{E:gamma}
	\bigg|\frac{1}{2\pi i} \int_{-\infty}^{(0+)} (\sqrt z)^{n-1} a^{-\sqrt z} \, e^{2zt} \, dz \bigg| =  
	O\bigg\{ \bigg(\frac{n}{2} \bigg)^{n/2} \bigg\},    
\end{equation}
\begin{equation}\label{E:gamma1}
\bigg|\frac{1}{2\pi i} \int_{-\infty}^{(0+)} \frac{(\sqrt z)^{n-1}}{\sqrt z + k} a^{-\sqrt z} \, e^{2zt} \, dz \bigg| =  O\bigg\{ \bigg(\frac{n}{2} \bigg)^{n/2} \bigg\},    
\end{equation}
where $n \to \infty$, and $a, t >0$ and $k = 1,2, \ldots$ are fixed.  
\end{Lem}
\begin{pf}
Estimate (\ref{E:gamma}) follows by integrating over the lower and the upper cut $\arg z = \mp \pi$, an obvious reparametrization, elementary estimates and Stirling's formula: 
\begin{equation}
\begin{align} \notag 
	&\bigg|\frac{1}{2\pi i} \int_{-\infty}^{(0+)} (\sqrt z)^{n-1} a^{-\sqrt z} \, e^{2zt} \, dz \bigg| < 
	2 \frac{1}{2\pi } \int^{\infty}_{0} (\sqrt x)^{n-1} |a^{\mp i \sqrt{x}}| \, e^{-2xt} \, dx = \\\notag 
	&2 \frac{1}{2\pi } \int^{\infty}_{0} (\sqrt x)^{n-1} \, e^{-2xt} \, dx <  2 \frac{1}{2\pi } \int^{\infty}_{0} x^{n/2-1} \, e^{-2xt} \, dx = \\  \notag
	&O\bigg\{\Gamma\bigg(\frac{n}{2}\bigg)\bigg\} = 
	O\bigg\{ \bigg(\frac{n}{2} \bigg)^{n/2} \bigg\}.  
\end{align}  
\end{equation}
Estimate (\ref{E:gamma1}) follows in a similar fashion.
\end{pf} 
\begin{Prop}\label{Th:premain} 
\begin{equation}\label{E:premain1}  
 	f_0(\lambda, t) = \sum_{m=0}^\infty d_m \, \frac{\sqrt \pi}{2\pi i} \int_{-\infty}^{(0+)} (\sqrt z)^{m-1} \, a^{-\sqrt z} \, e^{2zt} \, dz,    
\end{equation}
and for $n = 1,2, \ldots$
\begin{equation}\label{E:premain2}   
 	f_n(\lambda, t) = \sum_{k=0}^{[(n-1)/2]+1} a_k^{(n)} \sum_{m=0}^\infty d_m \, \frac{1}{2\pi i} 
	\int_{-\infty}^{(0+)}\frac{(\sqrt z)^{m}}{\sqrt z+k }\,a^{-\sqrt z} \, e^{2zt} \, dz,     
\end{equation}
where all the series are absolutely  convergent.
\end{Prop}
\begin{pf}
The proof follows by substituting the expansion 
\begin{equation}
	\frac{1}{ \Gamma(\sqrt z +1)} = \sum_{k=0}^{\infty} d_k \,  (\sqrt{z})^k,   
\end{equation}
in the integrals (\ref{E:f0x}) and (\ref{E:fnx}) and term by term integration.

The term by term integration is justifiable by absolute convergence and Fubini's theorem. The absolute convergence follows from estimates  (\ref{E:d}), (\ref{E:gamma}) and (\ref{E:gamma1}):

\begin{equation} \notag
\begin{align} \notag
 	&\bigg| \sum_{m=0}^\infty d_m \, \frac{1}{2\pi i} 
	\int_{-\infty}^{(0+)}\frac{(\sqrt z)^{m}}{\sqrt z+k }\,a^{-\sqrt z} \, e^{2zt} \, dz \bigg| <
	\sum_{m=0}^\infty |d_m| \, \bigg|\frac{1}{2\pi i} \int_{-\infty}^{(0+)}\frac{(\sqrt z)^{m}}{\sqrt z+k}\,a^{-\sqrt z} \, e^{2zt} \, dz\bigg|= 		\\ \notag
&= \sum_{m=0}^\infty O( m^{-(1-\varepsilon) m} ) \,  O\bigg\{ \bigg(\frac{m}{2} \bigg)^{m/2} \bigg\} < \infty, 
\end{align}    
\end{equation}
where  $k = 0, 1, \ldots$ and $ 0 <\varepsilon < 1/2$.
\end{pf} 
%
\begin{Prop}
Let $\alpha $ be a real parameter and $m = 0, 1, \ldots$.  \\
Then
\begin{equation}\label{E:He}
	\frac{1}{2\pi i } \int_{-\infty}^{(0+)}(\sqrt z)^{m-1}\,e^{- \alpha \sqrt z} \, e^{zt} \, dz  =
 	\frac{1}{\sqrt{\pi t}}\frac{1}{(2\sqrt t )^m} e^{-\alpha^2/(4t)} H_m \bigg( \frac{\alpha}{2\sqrt{t}} \bigg).	
\end{equation}
\begin{pf}
We compute initially the above integral in the particular case $m=0$, for which by the reparametrization $ z = w^2 $ and an elementary transform we obtain:
\begin{equation} \notag
	\begin{align} \notag
	&\frac{1}{2\pi i } \int_{-\infty}^{(0+)}(\sqrt z)^{-1}\,e^{- \alpha \sqrt z} \, e^{zt} \, dz =  
	\frac{1}{2\pi i } 2 \int_{c-i\infty}^{c+i\infty} e^{- \alpha  w} \, e^{ w^2 t} \, dw = \\ \notag 
	&= \frac{1}{2\pi i } 2 \int_{c-i\infty}^{c+i\infty} e^{- \alpha  w} \, e^{w^2 t} \, dw 
	= \frac{1}{\pi}  \int_{-\infty}^{\infty}  \cos{\alpha y} \,e^{-y^2 t} \, dy  = 
	\frac{1}{\sqrt{\pi t}} e^{-\alpha^2/(4t)}.
	\end{align}
\end{equation}
It is easy to see that the above relationship is differentiable any number of times on $\alpha $ under the sign of the integral, hence
\begin{equation} \notag
\begin{align}\label{E:Her}
	&\frac{(-1)^m}{2\pi i } \int_{-\infty}^{(0+)}(\sqrt z)^{m-1}\,e^{- \alpha \sqrt z} \, e^{zt} \, dz =  
	\frac{1}{\sqrt{\pi t}}\frac{\partial^m}{\partial \alpha^m}e^{-\alpha^2/(4t)} = \\ \notag
	&= \frac{1}{\sqrt{\pi t}} \frac{\partial^m}{\partial \alpha^m}e^{-(\alpha/(2\sqrt t))^2} =  
	\frac{1}{\sqrt{\pi t}}\frac{1}{(2\sqrt t )^m}\frac{d^m}{du^m}e^{-u^2} \bigg|_{u= \frac{\alpha}{2\sqrt t}} = \\ \notag
	&= \frac{1}{\sqrt{\pi t}}\frac{(-1)^m}{(2\sqrt t )^m} e^{-\alpha^2/(4t)} H_m \bigg( \frac{\alpha}{2\sqrt{t}} \bigg), 	
\end{align}	
\end{equation}
where the last equation has been written by virtue of the well--known relationship
\begin{equation} \notag
	H_m(x) = (-1)^m  e^{x^2} \frac{d^m}{du^m} (e^{-x^2});  
\end{equation}
(\ref{E:Her}) implies
\begin{equation} \notag
	\frac{1}{2\pi i } \int_{-\infty}^{(0+)}(\sqrt z)^{m-1}\,e^{- \alpha \sqrt z} \, e^{zt} \, dz =
	 \frac{1}{\sqrt{\pi t}}\frac{1}{(2\sqrt t )^m} e^{-\alpha^2/(4t)} H_m \bigg( \frac{\alpha}{2\sqrt{t}} \bigg), 	
\end{equation}
and this completes the proof. 
\end{pf}
\end{Prop}
\begin{Cor}\label{C:Her}
\begin{equation} \notag
	\frac{1}{2\pi i } \int_{-\infty}^{(0+)}(\sqrt z)^{m-1}\,a^{- \sqrt z} \, e^{2zt} \, dz =
	 \frac{1}{\sqrt{2\pi t}}\frac{1}{(2\sqrt{2t} )^m} e^{-\ln^2 a/(8t)} H_m \bigg( \frac{\ln a}{2\sqrt{2t}} \bigg). 	
\end{equation}
\end{Cor}
According to Bateman and Erd\'elyi (1954) {\bf 5.6} (12) 
we can write
\begin{equation}\label{E:Phi0}
 	\frac{1}{2\pi i} 
	\int_{-\infty}^{(0+)}\frac{e^{-\alpha \sqrt z}}{\sqrt z+ \beta }\, e^{zt} \, dz =
	\frac{e^{-\alpha^2/(4t)}}{\sqrt{\pi t}} - \beta e^{\alpha \beta + \beta^2t} \,
	erf\bigg(\frac{\alpha}{2\sqrt t}+ \beta \sqrt t \bigg),
\end{equation}
where $\alpha, \beta $ are real parameters. \\
Introduce the notation
\begin{equation}\label{E:Phimdef}
\begin{cases}
 	&\Phi_0(\alpha, \beta, t) :=  
	\frac{e^{-\alpha^2/(4t)}}{\sqrt{\pi t}} - \beta e^{\alpha \beta + \beta^2t} \,
	erf\bigg(\frac{\alpha}{2\sqrt t}+ \beta \sqrt t \bigg), \\
	&\Phi_m(\alpha, \beta, t) :=  (-1)^m \frac{\partial^m}{\partial \alpha^m}\Phi_0(\alpha, \beta, t),  
\end{cases}     
\end{equation}
where $m = 1,2, \ldots$. \\
{\em Remark}. 
Observe that repeated differentiation on $\alpha$ of the function  $\Phi_0(\alpha, \beta, t)$ results again in the error function 
of the same argument, $erf\bigg(\frac{\alpha}{2\sqrt t}+ \beta \sqrt t \bigg)$, which is multiplied by a certain exponential function with some exponential functions added. This follows from the definition
\begin{equation}\label{E:Phi0}
 	erf(x) = \frac{2}{\sqrt \pi} \int_0^x e^{-u^2} \,du
\end{equation}
and elementary rules of differentiation. The functions $\Phi_m(\alpha, \beta, t) $ are computable recursively by means of symbolic differentiation; for computations they are as good as elementary functions. 
\begin{Prop}
The integral in (\ref{E:Phi0}) as a function of the parameter $\alpha$ is differentiable under the sign of the integral any number of times, hence
\begin{equation}\label{E:Phim}
\begin{align}
	 &\Phi_m(\alpha, \beta, t) =  (-1)^m \frac{\partial^m}{\partial \alpha^m}\Phi_0(\alpha, \beta, t) = (-1)^m \frac{1}{2\pi i} \,  \frac{\partial^m}{\partial \alpha^m}
	\int_{-\infty}^{(0+)}\frac{ e^{-\alpha \sqrt z}}{\sqrt z+ \beta }\, e^{zt} \, dz = \\ \notag
	&= \frac{1}{2\pi i} \int_{-\infty}^{(0+)}\frac{(\sqrt z)^{m} e^{-\alpha \sqrt z}}{\sqrt z+ \beta }\, e^{zt} \, dz, 
\end{align}     
\end{equation}
where  $m = 1,2, \ldots$.
\end{Prop}
\begin{pf}
The last integral in (\ref{E:Phim}) absolutely converges for any $m$. The proof of  (\ref{E:Phim}) follows from (\ref{E:Phi0}) by induction on $m$ and by a standard application of bounded convergence theorem. 
\end{pf}
\begin{Cor}\label{C:erf}
\begin{equation} \notag
 \frac{1}{2\pi i} \int_{-\infty}^{(0+)}\frac{(\sqrt z)^{m} a^{- \sqrt z}}{\sqrt z+ k }\, e^{2zt} \, dz =
	\Phi_m(\ln a, k, 2t). 	
\end{equation}
\end{Cor}

The following theorem completes Laplace inversion; it is the main theorem of the present article.
\begin{Thm} \label{Th:main}
The density distribution $f(\lambda, t)$ of functional (\ref{E:fun}) is determined by the following 
absolutely convergent series:
\begin{equation} \label{E:fs}
f(\lambda, t) = \frac{1}{2\lambda \sqrt{2\lambda}} \,  e^{-1/(2\lambda)} e^{t/2}\sum_{n=0}^{\infty} \frac{1}{n!}\,
\, f_n(\lambda, t)\frac{1}{(2\lambda)^n},
\end{equation}
where
\begin{equation}\label{E:f0} 
	 f_0(\lambda, t) = \frac{1}{\sqrt{2t}} e^{-\ln^2{a}/(8t)} 
	\sum_{m=0}^\infty \frac{d_m}{(2\sqrt{2t})^m} \, H_m\bigg( \frac{\ln{a}}{2\sqrt{2t}} \bigg),
\end{equation}
and for $n = 1,2, \ldots$
\begin{equation}\label{E:fn} 
\begin{align} 
	f_n(\lambda, t) &= \\ \notag  
	&=\sum_{m=0}^\infty d_m \bigg\{ a_0^{(n)} \frac{e^{-\ln^2{a}/(8t)}}{\sqrt{2\pi t}\,(2\sqrt{2t})^m} 
	\, H_m\bigg( \frac{\ln{a}}{2\sqrt{2t}} \bigg) + \sum_{k=1}^{[(n-1)/2]+1} a_k^{(n)} \Phi_m(\ln{a}, k, 2t) \bigg \},   
	%
	%
\end{align} 
\end{equation}
where $a = 8\lambda e^{-2t}$, $H_m(x)$ are Hermite polynomials, $\Phi_m(\alpha, \beta, t)$
 are defined by (\ref{E:Phimdef}), the coefficients $d_m$ are defined by (\ref{E:dm1}) -- (\ref{E:dm2}), and the  coefficients $a_k^{(n)}$ are defined in Theorem \ref{Th:ank}.
\end{Thm}
\begin{pf}
The expansions above and their absolute convergence follow from Theorem \ref{Th:fser} and Proposition \ref{Th:premain} by substitutions indicated in Corollary \ref{C:Her} and Corollary \ref{C:erf}. \\ 
\end{pf}
{\em Remark}. Absolute convergence of these series can be also seen from the estimates of their remainders given in the next section, see 
Theorem \ref{Th:aserr} and Corollary \ref{C:2}. 
%
%
\section{Asymptotic Error Bounds} \label{S:err}
Estimates similar to those used in Section \ref{S:aux} readily lead to accurate error bounds that allow computations by series given in Theorem \ref{Th:main} with guaranteed accuracy. This matter, however, exceeds the limits of the present article.
In practice the required accuracy is usually achieved by other means of numerical analysis. For these reasons we content here with asymptotic estimates of error bounds. \\
The following theorem estimates the remainders of the series in 
(\ref{E:fs}) -- (\ref{E:fn}). \\
Let $R_N(\lambda,t)$,
$Q_M(\lambda,t)$ and $Q_M^{(k)}(\lambda,t)$,
where $k =  1, \ldots, [(n-1)/2]+1 $  be these remainders:
\begin{equation} \label{E:fsr} \notag
R_N(\lambda, t) = 
\sum_{n=N}^{\infty} \frac{1}{n!}\,
\, f_n(\lambda, t)\frac{1}{(2\lambda)^n},
\end{equation}
%
%
\begin{equation}\label{E:fnr}\notag 
 Q_M(\lambda, t) =   %
\sum_{m=M}^\infty \frac{d_m}{(2\sqrt{2t})^m} \, H_m\bigg( \frac{\ln{a}}{2\sqrt{2t}} \bigg) 
\end{equation}
\begin{equation}\label{E:fnkr} \notag
 Q_M^{(k)}(\lambda, t) = \sum_{m=M}^\infty d_m \Phi_m(\ln{a}, k, 2t).
\end{equation}
\begin{Thm} \label{Th:aserr}
 For fixed $\lambda, t> 0$    and $N, M \to \infty $ holds:    \\
(i)
\begin{equation} \label{E:fsr} \notag
| R_N(\lambda, t) |=  O \bigg(\frac{1}{(2\lambda)^N N!} \bigg), 
\end{equation}
(ii)
\begin{equation}\label{E:fnr} \notag 
 | Q_M(\lambda, t) | = O\big( M^{-(1-\varepsilon) M/2} \big),
\end{equation}
(iii)
\begin{equation}\label{E:fnr}\notag 
 | Q_M^{(k)}(\lambda, t) | = O\big( M^{-(1-\varepsilon) M/2} \big),
\end{equation}
where $k =  1, \ldots, [(n-1)/2]+1 $ is fixed.  
\end{Thm}
\begin{pf} Proof (i).
\begin{equation} \notag 
\begin{align} \notag
|R_N(\lambda, t)| &= 
	\bigg| \sum_{n=N}^{\infty} \frac{1}{n!} \, f_n(\lambda, t)\frac{1}{(2\lambda)^n} \bigg|
	 < \sum_{n=N}^{\infty} \frac{1}{n!}\,|f_n(\lambda, t)|\frac{1}{(2\lambda)^n},
\end{align} 
\end{equation}
where by (\ref{E:fn0}), (\ref{E:wn0}) and according to  (\ref{E:gest})

\begin{equation} \notag
\begin{align} \notag
	|f_n(\lambda, t)| &= \bigg| \frac{1}{2\pi i} \int_{c-i\infty}^{c+i\infty} w_n(\lambda, q) \, e^{qt}\,dq \bigg|=
	 \bigg|\frac{1}{2\pi i} \int_{c-i\infty}^{c+i\infty}  (2\lambda)^{-\sqrt {2q}/2} 
	\frac{\Gamma\big( \frac{\sqrt {2q}}{2}+n \big)}{\Gamma(\sqrt {2q}+n+1)} e^{qt}\, dq\bigg| < \\ \notag
	&<\frac{1}{2\pi } \int_{c-i\infty}^{c+i\infty}  |(2\lambda)^{-\sqrt {2q}/2}| 
	2\sqrt \pi \,\frac{1}{|\sqrt{2q}|}\big|2^{-\sqrt{2q}} \,\big| \,
	\bigg|\frac{1}{\Gamma(\sqrt{2q}/2 +1/2)} \bigg|.
	%
	%
 	 |e^{qt}| \, |dq| := \\ \notag
	&:= C(\lambda,t) < \infty, 
\end{align} 
\end{equation}
where $C(\lambda,t)$  does not depend on $n$. \\
This implies
\begin{equation} \notag
| R_N(\lambda, t) | < C(\lambda,t)\sum_{n=N}^{\infty} \frac{1}{n!} \frac{1}{(2\lambda)^n} = 
\frac{C(\lambda,t)}{(2\lambda)^N} \frac{1}{N!} O(1) =  O \bigg(\frac{1}{(2\lambda)^N N!} \bigg).  
\end{equation}
%
%
Proof (ii). 
By (\ref{E:premain1}) and estimates (\ref{E:d}) and  (\ref{E:gamma}) we have
\begin{equation}\label{E:f0r} 
	\begin{align} \notag
	 | Q_M(\lambda, t) | = 
%
%
\sum_{m=M}^\infty |d_m| \,\bigg|\frac{1}{2\pi i} 
	\int_{-\infty}^{(0+)} (\sqrt z)^{m-1/2}\,a^{-\sqrt z} \, e^{2zt} \, dz \bigg| = O\big( M^{-(1-\varepsilon) M/2} \big). 
	\end{align} 
\end{equation}
%
Proof (iii). In a similar fashion to Proof (ii) by (\ref{E:premain1}) and estimates (\ref{E:d}) and  (\ref{E:gamma}) we have
\begin{equation}\label{E:f0r} 
	\begin{align} \notag
	 | Q_M^{(k)}(\lambda, t) |&= \bigg| \sum_{m=1}^\infty d_m \, \frac{1}{2\pi i} 
\int_{-\infty}^{(0+)}\frac{(\sqrt z)^{m}}{\sqrt z+k }\,a^{-\sqrt z} \, e^{2zt} \, dz \bigg|   <  \\ \notag
	 &< \sum_{m=0}^\infty |d_m| \, \bigg|\frac{1}{2\pi i}
	\int_{-\infty}^{(0+)}\frac{(\sqrt z)^{m}}{\sqrt z+k }\,a^{-\sqrt z} \, e^{2zt} \, dz \bigg| =  \\ \notag 	
	&= O\big( M^{-(1-\varepsilon) M/2} \big). 
	\end{align} 
\end{equation}
\end{pf}
\begin{Cor} \label{C:2}
All the series in Theorem  \ref{Th:main} converge absolutely.
\end{Cor}
These error bounds show that for $\lambda$'s not too small the convergence rate of all the series in question is exceptionally good; for the small $\lambda$'s, however, an asymptotic expansion is available to fill this gap. This expansion will be discussed in a paper to follow.
%
\section{Computations} \label{S:ft}

The above series and error bounds provide effective means for computations of  
$f(\lambda,t)$, as illustrated on the example of the density distribution $f(\lambda,1)$ of the functional $\int_{0}^{1} e^{2 B(s)} ds$,
see Table 1. 
A graph of $f(\lambda,1)$ and further information on this distribution are available via {\em www.tolmatz.net}.

\begin{table}
\begin{tabular}{|c||c|c|c|c|c|c|c|}
\hline
	$\lambda$ &$f(\lambda,1)$ &$\lambda$ &$f(\lambda,1)$&$\lambda$&$f(\lambda,1)$&$\lambda$ &$f(\lambda,1)$ \\  
 \hline 
 0.01&0.0000000000&0.51&0.5815736451&1.01& 0.3469311307&1.51&0.2149882566   \\   
 0.02&0.0000000009&0.52&0.5768409044& 1.02&0.3433418909&1.52&0.2131233054    \\ 
 0.03&0.0000025112&0.53&0.5719894905  &   1.03&0.3398003064&1.53&0.2112812233 \\     
 0.04&0.0001196142&0.54&0.5670368103&1.04   & 0.3363058318&1.54&0.2094616620 \\     
 0.05&0.0011483687&0.55&0.5619986842&1.05& 0.3328579122&1.55&0.2076642790     \\ 
         0.06&0.0049930908&0.56&0.5568894834&1.06&
    0.3294559852&1.56&0.2058887374      \\
         0.07&0.0138753955&0.57&0.5517222559&1.07&
    0.3260994820&1.57&0.2041347063      \\
         0.08&0.0292395656&0.58&0.5465088421&1.08&
    0.3227878293&1.58&0.2024018600      \\
         0.09&0.0513458853&0.59&0.5412599790&1.09&
    0.3195204500&1.59&0.2006898783      \\
         0.10&0.0794780720&0.60&0.5359853972&1.10&
    0.3162967645&1.60&0.1989984467      \\
         0.11&0.1123591601&0.61&0.5306939077&1.11&
    0.3131161920&1.61&0.1973272557      \\
         0.12&0.1485212098&0.62&0.5253934820&1.12&
    0.3099781508&1.62&0.1956760010      \\
         0.13&0.1865502642&0.63&0.5200913256&1.13&
    0.3068820597&1.63&0.1940443838      \\
         0.14&0.2252159442&0.64&0.5147939437&1.14&
    0.3038273385&1.64&0.1924321099      \\
         0.15&0.2635203296&0.65&0.5095072027&1.15&
    0.3008134089&1.65&0.1908388905      \\
         0.16&0.3006989269&0.66&0.5042363855&1.16&
    0.2978396947&1.66&0.1892644415      \\
         0.17&0.3361973070&0.67&0.4989862419&1.17&
    0.2949056229&1.67&0.1877084836      \\
         0.18&0.3696381349&0.68&0.4937610351&1.18&
    0.2920106239&1.68&0.1861707424      \\
         0.19&0.4007868533&0.69&0.4885645839&1.19&
    0.2891541317&1.69&0.1846509478      \\
         0.20&0.4295201362&0.70&0.4834003011&1.20&0.2863355851&
    1.70&0.1831488348      \\
         0.21&0.4557987888&0.71&0.4782712290&1.21&
    0.2835544271&1.71&0.1816641425      \\
         0.22&0.4796454414&0.72&0.4731800715&1.22&
    0.2808101060&1.72&0.1801966146      \\
         0.23&0.5011267281&0.73&0.4681292234&1.23&
    0.2781020752&1.73&0.1787459991      \\
         0.24&0.5203393667&0.74&0.4631207976&1.24&
    0.2754297938&1.74&0.1773120484      \\
         0.25&0.5373994870&0.75&0.4581566492&1.25&  
0.2727927265&1.75&0.1758945189       \\ \hline                 
\end{tabular}
\end{table}

\begin{table}
\begin{tabular}{|c||c|c|c|c|c|c|c|}
\hline
	$\lambda$ &$f(\lambda,1)$ &$\lambda$ &$f(\lambda,1)$&$\lambda$&$f(\lambda,1)$&$\lambda$ &$f(\lambda,1)$ \\  
 \hline 
         0.26&0.5524345856&0.76&0.4532383988&1.26&
    0.2701903441&1.76&0.1744931712      \\
         0.27&0.5655775645&0.77&0.4483674521&1.27&
    0.2676221234&1.77&0.1731077702      \\
         0.28&0.5769623929&0.78&0.4435450197&1.28&
    0.2650875475&1.78&0.1717380844      \\
         0.29&0.5867210214&0.79&0.4387721336&1.29&
    0.2625861059&1.79&0.1703838864      \\
         0.30&0.5949812497&0.80&0.4340496636&1.30&0.2601172946&
    1.80&0.1690449526      \\
         0.31&0.6018653139&0.81&0.4293783310&1.31&
    0.2576806161&1.81&0.1677210634      \\
         0.32&0.6074890101&0.82&0.4247587227&1.32&
    0.2552755797&1.82&0.1664120026      \\
         0.33&0.6119612147&0.83&0.4201913026&1.33&
    0.2529017013&1.83&0.1651175578      \\
         0.34&0.6153836919&0.84&0.4156764232&1.34&
    0.2505585033&1.84&0.1638375201      \\
         0.35&0.6178511084&0.85&0.4112143355&1.35&
    0.2482455154&1.85&0.1625716841      \\
         0.36&0.6194511911&0.86&0.4068051983&1.36&
    0.2459622734&1.86&0.1613198481      \\
         0.37&0.6202649828&0.87&0.4024490870&1.37&
    0.2437083205&1.87&0.1600818135      \\
         0.38&0.6203671604&0.88&0.3981460013&1.38&
    0.2414832062&1.88&0.1588573850      \\
         0.39&0.6198263907&0.89&0.3938958721&1.39&
    0.2392864869&1.89&0.1576463709      \\
         0.40&0.6187057040&0.90&0.3896985686&1.40&0.2371177259&
    1.90&0.1564485823      \\
         0.41&0.6170628736&0.91&0.3855539041&1.41&
    0.2349764928&1.91&0.1552638337      \\
         0.42&0.6149507900&0.92&0.3814616412&1.42&
    0.2328623643&1.92&0.1540919427      \\
         0.43&0.6124178242&0.93&0.3774214976&1.43&
    0.2307749233&1.93&0.1529327297      \\
         0.44&0.6095081762&0.94&0.3734331500&1.44&
    0.2287137596&1.94&0.1517860184      \\
         0.45&0.6062622048&0.95&0.3694962388&1.45&
    0.2266784693&1.95&0.1506516352      \\
         0.46&0.6027167376&0.96&0.3656103720&1.46&
    0.2246686549&1.96&0.1495294094      \\
         0.47&0.5989053617&0.97&0.3617751285&1.47&
    0.2226839255&1.97&0.1484191731      \\
         0.48&0.5948586930&0.98&0.3579900616&1.48&
    0.2207238963&1.98&0.1473207614      \\
         0.49&0.5906046267&0.99&0.3542547020&1.49&
    0.2187881889&1.99&0.1462340117      \\
         0.50&0.5861685681&1.00&0.3505685606&1.50&0.2168764311&
    2.00&0.1451587645   \\ \hline                
\end{tabular}
\caption{The density function $f(\lambda,1)$.}
\end{table}

\end{document}